\documentclass[12pt]{amsart}
\usepackage{amsmath,amscd,amssymb,amsfonts}
\def\ver{May 27, 2007, v.1}

\def\ssbull{\raise.2ex\hbox{${\scriptscriptstyle\bullet}$}}
\def\scirc{\,\raise.2ex\hbox{${\scriptstyle\circ}$}\,}
\def\mcup{\hbox{$\bigcup$}}

\def\msum{\hbox{$\sum$}}
\def\mopl{\hbox{$\bigoplus$}}
\def\fa{{\mathfrak a}}
\def\fm{{\mathfrak m}}
\def\be{{\mathbf e}}
\def\bC{{\mathbb C}}
\def\bD{{\mathbb D}}
\def\bN{{\mathbb N}}
\def\bQ{{\mathbb Q}}
\def\bR{{\mathbb R}}
\def\bZ{{\mathbb Z}}
\def\boR{{\mathbf R}}
\def\cD{{\mathcal D}}
\def\cG{{\mathcal G}}
\def\cH{{\mathcal H}}
\def\cI{{\mathcal I}}
\def\cJ{{\mathcal J}}
\def\cM{{\mathcal M}}
\def\cO{{\mathcal O}}
\def\cX{{\mathcal X}}
\def\ocA{\bar{\mathcal A}}
\def\oA{\bar{A}}
\def\oB{\bar{B}}
\def\oG{\overline{G}}
\def\oI{\bar{I}}
\def\tb{\tilde{b}}
\def\ti{\tilde{i}}
\def\tE{\widetilde{E}}
\def\tH{\tilde{H}}
\def\tY{\widetilde{Y}}
\def\tZ{\widetilde{Z}}
\def\hsig{\widehat{\sigma}}
\def\red{{\rm red}}
\def\rank{{\rm rank}}
\def\Gr{{\rm Gr}}
\def\JC{{\rm JC}}
\def\Spec{{\rm Spec}}
\def\cSpec{{{\mathcal S}pec}}
\def\cSpecan{{{\mathcal S}pecan}}
\def\Int{{\rm Int}\,}
\def\Sing{{\rm Sing}\,}
\def\supp{{\rm supp}}
\def\Sp{{\rm Sp}}
\def\hSp{\widehat{\rm S}{\rm p}}
\def\Ker{{\rm Ker}}
\def\mp{\,\widetilde{+}\,}

\begin{document}
\title{Spectrum and multiplier ideals of arbitrary subvarieties}
\author{Alexandru Dimca }
\address{Laboratoire J.A. Dieudonn\'e, UMR du CNRS 6621,
Universit\'e de Nice-Sophia Antipolis, Parc Valrose,
06108 Nice Cedex 02, FRANCE.}
\email{Alexandru.DIMCA@unice.fr}
\author{Philippe Maisonobe}
\address{Laboratoire J.A. Dieudonn\'e, UMR du CNRS 6621,
Universit\'e de Nice-Sophia Antipolis, Parc Valrose,
06108 Nice Cedex 02, FRANCE.}
\email{Philippe.MAISONOBE@unice.fr}
\author{Morihiko Saito }
\address{RIMS Kyoto University, Kyoto 606--8502 JAPAN}
\email{ msaito@kurims.kyoto-u.ac.jp}
\date{\ver}
\begin{abstract}
We introduce a spectrum for arbitrary varieties.
This generalizes the definition by Steenbrink for hypersurfaces.
In the isolated complete intersection singularity case,
it coincides with the one given by Ebeling and Steenbrink
except for the coefficients of integral exponents.
We show a relation to the graded pieces of the multiplier ideals
by using a relation to the filtration $V$ of Kashiwara and Malgrange.
This implies a partial generalization of a theorem of Budur
in the hypersurface case.
The point is to consider the direct sum of the graded pieces of
the multiplier ideals as a module over the algebra defining
the normal cone of the subvariety.
We also give a combinatorial description in the case of monomial
ideals.
\end{abstract}
\maketitle

\bigskip\centerline{\bf Introduction}

\bigskip\noindent
In [27], [28], Steenbrink introduced the spectrum for hypersurface
singularities.
Its relations with $b$-function and multiplier ideals have been studied
in [3], [6], [24], [25], etc.
The multiplier ideals were originally defined for any subvariety of a
smooth variety (see e.g. [16]), and the $b$-function for an arbitrary
variety has been defined in [4].
In this paper we introduce the spectrum for an arbitrary variety
generalizing Steenbrink's definition [27], [28].

Let $X$ be a closed subvariety of a smooth complex algebraic variety
or a complex manifold $Y$.
Let $(N_XY)_x$ be the fiber of the normal cone $N_XY\to X$ over
$x\in X$.
For each irreducible component $\Lambda$ of $(N_XY)_x$, set
$n_{\Lambda}=\dim Y-\dim\Lambda$.
We have the nonreduced spectrum and the (reduced) spectrum
$$
\hSp(X,\Lambda)=\msum_{\alpha>0}\,m_{\Lambda,\alpha}t^{\alpha}\in
\bZ[t^{1/e}],\quad
\Sp(X,\Lambda)=\hSp(X,\Lambda)-(-1)^{n_{\Lambda}}t^{n_{\Lambda}+1},
$$
where $e\in\bZ_{>0}$, see (1.2) for the definition of
$m_{\Lambda,\alpha}$.
Note that $\hSp(X,\Lambda)$, $\Sp(X,\Lambda)$ are essentially
independent of $Y$ as a corollary of a product formula,
see Cor.~(3.3) and (3.4).
In case $(N_XY)_x$ is irreducible, set $\hSp(X,x)=\hSp(X,\Lambda)$
for $\Lambda=(N_XY)_x$, and similarly for $\Sp(X,x)$, $m_{x,\alpha}$.
This generalizes Steenbrink's definition in the hypersurface case
where $N_XY$ is a line bundle over $X$.
We use $\Sp(X,x)$ mainly in this case.
The difference between $\hSp(X,\Lambda)$ and $\Sp(X,\Lambda)$
comes from the one between cohomology and reduced cohomology.

In the isolated complete intersection singularity case, $N_XY$ is a
vector bundle over $X$ (in particular, $(N_XY)_x$ is irreducible),
and the spectrum is associated to the mixed Hodge structure on the
Milnor cohomology where the action of the monodromy is given by
choosing a sufficiently general line passing through the origin in
the base space of the Milnor fibration,
see also Remark~(1.3)(i).
In this case Ebeling and Steenbrink [10] defined the spectrum
in a different way where they consider also the contribution
of the Milnor cohomology associated to the singularity of the
total space of a generic 1-parameter smoothing of $X$
so that the symmetry and the semicontinuity hold.
Their spectrum differs from ours in general ([10], [19], [29]),
but they coincide for $m_{x,\alpha}$ with $\alpha\notin\bZ$,
see Remark~(1.3)(iv).
So we can apply Theorem~1 below to their spectrum except for
the case $\beta=1$.

Let $\cJ(Y,\alpha X)$ denote the multiplier ideals of $X\subset Y$ for
$\alpha\in\bQ_{>0}$.
They can be defined by the local integrability of
$$
|g|^2/(\msum_i\,|f_i|^2)^{\alpha}\,\,\,\,\text{for}\,\,\, g\in\cO_Y,
$$
where $f_1,\dots,f_r$ are local generators of the ideal of $X$,
see [16].
They are closely related to the filtration $V$ of Kashiwara [14] and
Malgrange [18], see [4].
Set
$$
\cG(Y,\alpha X)=\cJ(Y,(\alpha-\varepsilon) X)/\cJ(Y,\alpha X)\,\,\,
\text{for $0<\varepsilon\ll 1$}.
$$
If $\cG(Y,\alpha X)_x\ne 0$, then $\alpha$ is called a jumping
coefficient of $X\subset Y$ at $x$.

Let $\beta\in(0,1]\cap\bQ$, and $\cI_X$ be the ideal sheaf of
$X\subset Y$.
Define
$$
\cM(\beta)=\mopl_{i\in\bN}\,\cG(Y,(\beta+i)X),\quad
\ocA=\mopl_{i\in\bN}\,\cI_X^i/\cI_X^{i+1}.
$$
Then $\cM(\beta)$ is a graded $\ocA$-module since
$\cI_X\cJ(Y,\alpha X)\subset\cJ(Y,(\alpha+1)X)$.
Set
$$
\tZ_{\beta}=\supp_{\ocA}\,\cM(\beta)\subset N_XY=\cSpec_X\ocA,
$$
i.e. $\tZ_{\beta}$ is the support of the associated sheaf on
$N_XY$.
Here $\cSpec_X\ocA$ is replaced by $\cSpecan_X\ocA$ in the analytic
case.
For an irreducible component $E$ of $N_XY$, let
$E^0$ be the complement in $E$ of the intersection of $E$ with
the union of the other irreducible components of $N_XY$.
Let $\tZ_{\beta,E}$ be the closure of $\tZ_{\beta}\cap E^0$.
Let $Z_{\beta,E}$ be the image of $\tZ_{\beta,E}$ in $X$.
Set $d_{\beta,E}=\dim Z_{\beta,E}$.

For $x\in X$, let $\fm_{Y,x}$ be the maximal ideal of $\cO_{Y,x}$,
and set
$$
\cM(\beta,x)=\cM(\beta)/\fm_{Y,x}\cM(\beta),\quad
\ocA(x)=\ocA/\fm_{Y,x}\ocA.
$$
Then $\cM(\beta,x)$ is a graded $\ocA(x)$-module.
For each irreducible component $\Lambda$ of $(N_XY)_x=\Spec\,\ocA(x)$,
let $\mu_{\Lambda,\beta}$ be the multiplicity of $\cM(\beta,x)$ at
the generic point of $\Lambda$, i.e.
$$
\mu_{\Lambda,\beta}=\dim_{K(\Lambda)}\cM(\beta,x)\otimes_{\ocA(x)}K
(\Lambda),
$$
where $K(\Lambda)$ is the function field of $\Lambda$ which is a
localization of $\ocA(x)_{\red}$.

\medskip\noindent
{\bf Theorem~1.} {\it Let $\beta\in(0,1]\cap\bQ$. Then

\smallskip\noindent
{\rm (i)} We have in general
$$
0\le m_{\Lambda,\beta}\le\mu_{\Lambda,\beta}.
$$
In particular, $m_{\Lambda,\beta}=0$ if
$x\notin\supp\,\cM(\beta)\subset X$.

\smallskip\noindent
{\rm (ii)} If $\tZ_{\beta}$ is contained in $\Lambda$ on a
neighborhood of $\xi$ in $N_XY$, then
$$
m_{\Lambda,\beta}=\mu_{\Lambda,\beta}.
$$

\smallskip\noindent
{\rm (iii)}
If $x$ is a sufficiently general point of $Z_{\beta,E}$,
then for any irreducible component $\Lambda$ of $(N_XY)_x$
contained in $E$, we have $m_{\Lambda,\beta+i}=0$  for any integer
$i<d_{\beta,E}$, and
}
$$
m_{\Lambda,\beta+d_{\beta,E}}=(-1)^{d_{\beta,E}}\mu_{\Lambda,\beta}.
$$

\medskip
Here a sufficiently general point means that it belongs to a
(sufficiently small) non-empty Zariski-open subset of $Z_{\beta,E}$.
This gives a partial generalization of a theorem of Budur~[3]
in the hypersurface case where $\cM(\beta)$ is a free $\ocA$-module
of rank 1 over $\cG(Y,\beta X)$.
Theorem~1 was found during an attempt to extend some assertions about
the spectrum in [9].
The following is a generalization of Cor~1.5 and 1.6 in loc.~cit.

\medskip\noindent
{\bf Theorem~2.} {\it Let $\{S_i\}$, $\{S'_j\}$ be Whitney
stratifications of $X$, $N_XY$ such that the restriction of the
projection $N_XY\to X$ to each $S'_j$ is a smooth morphism to
some $S_i$, and the restriction of $\cH^{\ssbull}\Sp_X\bC_Y$ to
$S'_j$ are local systems where $\Sp_X\bC_Y$ denotes the Verdier
specialization {\rm [31]}. Then

\smallskip\noindent
{\rm (i)} The spectrum $\hSp(X,\Lambda_x)$ remains constant if
$\{\Lambda_x\}_{x\in U}$ is a locally trivial family of
irreducible components of $(N_XY)_x$ for $x\in U\subset S_i$
where $U$ is an analytically open subset of $S_i$.

\smallskip\noindent
{\rm (ii)} If $T$ is a transversal slice to $S_i$ in $Y$ such that
$S_i\cap T=\{x\}$, then we have
$$
\hSp(X,\Lambda)=(-t)^{d_i}\,\hSp(X\cap T,\Lambda),
$$
for an irreducible component $\Lambda$ of
$(N_{X\cap T}T)_x=(N_XY)_x$ where $d_i=\dim S_i$.
}

\medskip
Here a locally trivially family $\{\Lambda_x\}_{x\in U}$ in (i)
means a local section of the sheaf defined by the sets of the
irreducible components of $(N_XY)_x$ for $x\in S_i$, which is
a locally constant sheaf.
Theorem~2 implies that Theorem~1(iii) is reduced to (ii)
(where $d_{\beta,E}=0$) by restricting to a sufficiently
general member of a family of transversal slices.
Note that for $\beta\notin\bZ$, we have
$Z_{\beta,E}\subset\supp\,\cM(\beta)\subset\Sing X$,
and hence $d_{\beta,E}=0$ in the isolated singularity case.

Let $f=(f_1,\dots,f_r)$ be a set of local generators of $\cI_X$, and
$b_f(s)$ be the $b$-function of $f$ in the sense of [4] (see also
[11], [20]).
This is independent of the choice of $f$ and $r$,
but depends on the choice of $Y$.
By [4] and Theorem~1 we get

\medskip\noindent
{\bf Theorem~3.} {\it
Let $\alpha\in\bQ_{>0}$, and assume $m_{\Lambda,\alpha}\ne 0$
for some irreducible component $\Lambda$ of $(N_XY)_x$.
Then $\alpha +i$ is a root of $b_f(-s)$ for some $i\in\bZ$.
If furthermore $\alpha<1$, then there is a nonnegative integer
$j_0$ such that $\alpha+j$ is a jumping coefficient
of $X\subset Y$ at $x$ for any integer $j\ge j_0$.
}

\medskip
It is not easy to determine $i$ and $j_0$ in Theorem~3
(see Example~(4.6) below) except for $j_0$ in
the hypersurface case
(here $j_0=0$ since $\cG(Y,\alpha X)_x=\cG(Y,(\alpha+j)X)_x$ for
$\alpha>0$ and $j\in\bN$).
In the monomial ideal case, $j_0$ is bounded by $\dim Y-1$,
see Remark~(4.8) below.

In the monomial ideal case, there is a combinatorial description
for the jumping coefficients [12] and for the roots of the
$b$-function [5].
We give here one for the spectrum, see Th.~(4.4).

In Section~1, we review the specializations and define the spectrum.
In Section~2, we prove Theorems~1--3.
In Section~3, we show a product formula which implies the
well-definedness of the spectrum.
In Section~4, we treat the monomial ideal case, and prove Th.~(4.4).

\bigskip\bigskip
\centerline{\bf 1. Spectrum}

\bigskip
In this section we review the specializations and define the spectrum.

\medskip\noindent
{\bf 1.1.~Specialization.}
Let $Y$ be a smooth complex algebraic variety or a complex
manifold, and $X$ be a subvariety or a analytic subspace of $Y$.
We do not assume $X$ reduced nor irreducible.
Let $N_XY$ denote the normal cone.
Let $f=(f_1,\dots,f_r)$ be a set of local generators of the ideal
$\cI_X$ of $X\subset Y$. We denote the graph embedding by
$$
i_f:Y\to \tY:=Y\times\bC^r.
$$
Let $z_1,\dots,z_r$ be the coordinates of $\bC^r$.
We have the canonical surjection
$$
\ocA':=\cO_Y[z_1,\dots,z_r]\to\ocA:=\mopl_{i\in\bN}\,
\cI_X^i/\cI_X^{i+1},
$$
sending $z_i$ to $f_i$.
This implies the inclusion $N_XY\subset \tY$ such that
the projection $\pi':\tY\to Y$ induces $\pi:N_XY\to X$.
Let $\partial_i=\partial/\partial z_i$, and define
$$
(M,F):=(i_f)_*(\cO_Y,F)=(\cO_Y[\partial_1,\dots,\partial_r],F).
$$
This is a direct image as a filtered $\cD$-module.
Setting $\partial^{\nu}=\prod_i \partial_i^{\nu_i}$, we have
$$
F_pM=\msum_{|\nu|\le p+\dim Y}\cO_Y\otimes\partial^{\nu}.
\leqno(1.1.1)
$$

Let $V$ denote the filtration of Kashiwara [14] and Malgrange [18]
on $M$ along $Y\times\{0\}$ indexed by $\bQ$.
The specialization of $(\cO_Y,F)$ along $X$ is
defined by
$$
\Sp_X(\cO_Y,F)=
\pi'{}^{-1}(\mopl_{\alpha\in\bQ}\,\Gr_V^{\alpha}(M,F))
\otimes_{\pi'{}^{-1}\ocA'}\cO_{\tY}.
$$
By [4], [9] this is compatible with
the the definition of specialization [31] in the category of
perverse sheaves [1] or mixed Hodge modules [23]
$$
\Sp_X\bQ_Y[\dim Y]=\psi_t\boR j_*\bQ_{Y\times\bC^*}[\dim Y],
\leqno(1.1.2)
$$
where
$$
j:Y\times\bC^*=\cSpec_Y(\cO_Y[t,t^{-1}])
\to\cSpec_Y(\mopl_{i\in\bZ}\,\cI_X^{-i}\otimes t^i),
$$
denotes the inclusion to the deformation to the normal cone, see [31].
Here $\cI_X^{-i}=\cO_Y$ for $i\ge 0$, and $\cSpec_Y$ is replaced by
$\cSpecan_Y$ in the analytic case.
Note that the action of the semisimple part $T_s$ of the monodromy
associated to the nearby cycle functor $\psi_t$ corresponds to the
multiplication by $\exp(-2\pi i\alpha)$ on $\Gr_V^{\alpha}(M,F)$.
These follow from the fact (see [4], [9]) that the filtration $V$
corresponding to the functor $\psi_t$ is essentially given by
$$
\mopl_{i\in\bZ}\,V^{\alpha-i}M\otimes t^i.
$$

By [4] (see also [9] for the analytic case)
we have the isomorphisms
$$
F_{-\dim Y}\Gr_V^{\alpha}M=\cG(Y,\alpha X)\quad\text{for}\,\,\,
\alpha\in\bQ.
\leqno(1.1.3)
$$
For $\cM(\beta)$ as in the introduction, this implies
$$
\aligned
F_{-\dim Y}\Sp_X(\cO_Y)
&=\pi'{}^{-1}(\mopl_{0<\beta\le 1}\,\cM(\beta))
\otimes_{\pi'^{-1}\ocA'}\cO_{\tY}\\
&=\pi'{}^{-1}(\mopl_{0<\beta\le 1}\,\cM(\beta))
\otimes_{\pi^{-1}\ocA}\cO_{N_XY}.
\endaligned
$$
since $\cO_{N_XY}=\pi^{-1}\ocA\otimes_{\pi'^{-1}\ocA'}\cO_{\tY}$
for the last isomorphism.
These imply that $\cM(\beta)$ is locally finitely generated over
$\ocA'$ or $\ocA$.

\medskip\noindent
{\bf 1.2.~Definition.}
With the above notation, set
$$
(N_XY)_x=\pi^{-1}(x)\,\,\,\,\text{for}\,\,\,x\in X.
$$
Let $\Lambda$ be an irreducible component of $(N_XY)_x$, and
take a sufficiently general point $\xi$ of $\Lambda$.
Let $i_{\xi}:\{\xi\}\to N_XY$ denote the inclusion
morphism.
Set
$$
d_Y=\dim Y,\quad d_{\Lambda}=\dim\Lambda,\quad
n_{\Lambda}=d_Y-d_{\Lambda}.
$$
We denote the pull-back of $\Sp_X(\cO_Y,F)[-d_Y]$ as a complex
of filtered $\cD$-modules by
$$
i_{\xi}^*\Sp_X(\cO_Y,F)[-d_Y],
$$
see Remark~(1.3)(ii) below.
This corresponds to $i_{\xi}^*\Sp_X\bC_Y$,
and underlies a complex of mixed Hodge module on $\{\xi\}$,
which is identified with a complex of mixed Hodge structures [7],
see [23].
Combined with the action of  $T_s$, it defines the nonreduced
spectrum as in [27], [28]:
$$
\hSp(X,\Lambda)=\msum_{\alpha}\,m_{\Lambda,\alpha}t^{\alpha},
$$
where
$$
\aligned
&m_{\Lambda,\alpha}=\msum_j\,(-1)^j\dim\Gr_F^pH^{j+n_{\Lambda}}
(i_{\xi}^*\Sp_X\cO_Y[-d_Y])_{\lambda},\\
&\quad\text{with}\quad
p=[n_{\Lambda}+1-\alpha],\,\,\lambda=\exp(-2\pi i\alpha).
\endaligned
\leqno(1.2.1)
$$
Here $H^{\ssbull}(i_{\xi}^*\Sp_X\cO_Y[-d_Y])_{\lambda}$ denotes the
$\lambda$-eigenspace of the cohomology group by the action of $T_s$.
Set
$$
\Sp(X,\Lambda)=\hSp(X,\Lambda)-(-1)^{n_{\Lambda}}t^{n_{\Lambda}+1}.
$$
If $(N_XY)_x$ is irreducible, set for $\Lambda=(N_XY)_x$
$$
\hSp(X,x)=\hSp(X,\Lambda),\quad \Sp(X,x)=\Sp(X,\Lambda),\quad
m_{x,\alpha}=m_{\Lambda,\alpha}.
$$

\medskip\noindent
{\bf 1.3.~Remarks.}~(i) This generalizes Steenbrink's definition in
the hypersurface case ([27], [28]).
Assume $X$ is a hypersurface or an isolated complete intersection
singularity.
Then $N_XY$ is a line bundle or a vector bundle over $X$,
and the Milnor cohomology is given by
$$
H^{\ssbull}(i_{\xi}^*\Sp_X(\cO_Y,F)[-d_Y])\,\,\,\text{for}\,\,\,
\xi\in(N_XY)_x\,\,\,\text{sufficiently general},
$$
where the pull-back $i_{\xi}^*$ is explained in Remark~(ii) below.
In the isolated complete intersection singularity case,
the action of the monodromy is given by taking a sufficiently
general line $C$ passing through the origin in the base space of the
Milnor fibration and  corresponding to $\xi\in\Lambda=(N_XY)_{x}$.
Indeed, if the complete intersection $X$ is defined by
$$
f=(f_1,\dots,f_r):Y\to\bC^r,
$$
then $f$ induces the projection
$$
N_XY=X\times\bC^r\to N_{\{0\}}\bC^r=\bC^r,
$$
and the inverse image $Z$ of a sufficiently general line
$C\subset\bC^r$ by $f$ is a 1-parameter smoothing of $X$.
Moreover, its Milnor cohomology is isomorphic to the cohomology
of $(\Sp_X\bC_Y)_{\xi}$ with $\bC\,\xi\subset(N_XY)_x$
corresponding to $C$ by the above projection,
and the restriction of $\Sp_X\bC_Y$ to
$N_XZ\subset N_XY$ is isomorphic to $\Sp_X\bC_{Z}$ on a
neighborhood of $\xi\in(N_XZ)_x$ if $\xi$ is sufficiently general.
(Note also that we can replace $Y$ with the total space of the
miniversal deformation of $f$ by Cor.~(3.4).)

\medskip
(ii) Under a closed embedding of smooth complex varieties or complex
manifolds $i:X\to Y$, the pull-back of a complex of filtered
$\cD$-modules $(M,F)$ underlying a complex of mixed Hodge modules
is locally defined as follows.
By factorizing $i$ locally on $Y$, we may assume $X$ is defined by
$y_1=0$ where $y_1,\dots,y_m$ are local coordinates of $Y$.
Let $V$ be the filtration of Kashiwara [14] and Malgrange [18]
along $X$.
Then the pull-back $i^*(M,F)$ as a complex of filtered $\cD$-modules
is defined to be the mapping cone of
$$
\partial/\partial y_1:\Gr_V^1(M,F[1])\to\Gr_V^0(M,F),
\leqno(1.3.1)
$$
where $(F[1])_p=F_{p-1}$.
This is the same as in the case of right $\cD$-modules, and
the transformation between the corresponding left and right
$\cD$-modules is done without shifting the filtration in this paper.
(Note that the usual definition of pull-back as in [2] does not
work for the Hodge filtration.)

\medskip
(iii) Define
$$
p(M,F)=\min\{p\,|\,F_pM\ne 0\}.
\leqno(1.3.2)
$$
By Remark~(ii), $p(M,F)$ does not decrease under the cohomological
pull-back functor $\cH^ki*$ for filtered $\cD$-modules by a closed
embedding $i$, and increases by the codimension under
a non-characteristic restriction to a transversal slice
(since $\Gr_V^0M=0$ in the non-characteristic case).

Under the pull-back by a smooth morphism,
$p(M,F)$ decreases by the relative dimension
where the complex is also shifted.
This is compatible with the definition of $F$ on $\cO_Y$,
and we have
$$
p(\cO_Y,F)=p(\Sp_X\cO_Y,F)=-\dim Y.
\leqno(1.3.3)
$$
Then we get by Remark~(ii)
$$
m_{\Lambda,\alpha}=0\,\,\,\text{if}\,\,\,\alpha\le 0,\quad
m_{\Lambda,\beta}\ge 0\,\,\,\text{if}\,\,\,\beta\in(0,1].
\leqno(1.3.4)
$$

\medskip
(iv) The normalization of spectrum (1.2.1) is dual of the one used by
Ebeling and Steenbrink [10].
Indeed, for a Hodge structure $H$ with an automorphism $T$ of finite
order, they use
$$
\aligned
&\Sp'(H,T)=\msum_{\alpha}\,m_{H,\alpha}t^{\alpha}\quad
\text{with}\quad m_{H,\alpha}=\dim_{\bC}\Gr^p_FH_{\bC,\lambda},\\
&\qquad\text{where}\quad p=[\alpha],\,\,\,\lambda=\exp(2\pi i\alpha).
\endaligned
$$
Here $H_{\bC,\lambda}=\Ker(T-\lambda)\subset H_{\bC}$.
This definition is somewhat dual of (1.2.1).
In the case of isolated hypersurface singularities, their
definition of spectrum coincides with ours by the symmetry of the
spectrum which follows from the self-duality of the mixed Hodge
structure on the Milnor cohomology [27].
In the case of isolated complete intersection singularities,
they apply the above definition to the mixed Hodge structure
$$
\varphi_f\psi_g \bQ_{\cX}[d_X],
$$
where $f:X'\to\bC$ is a generic 1-parameter smoothing of $X$,
$g:\cX\to\bC$ is a generic 1-parameter smoothing of $X'$,
$T$ is the semisimple part of the monodromy associated to
$\varphi_f$, and $d_X=\dim X$, see [10] for details.
Then the symmetry of their spectrum follows from the self-duality
of $\varphi_f\psi_g \bQ_{\cX}[d_X]$ in [23], 2.6.2.

Denoting the Milnor fibers of $f,g$ by $F_f,F_g$,
we have a short exact sequence
$$
0\to\tH^{d_X}(F_f,\bC)\to\varphi_f\psi_g \bQ_{\cX}[d_X]\to
H^{d_X+1}(F_g,\bC)\to 0,
$$
since $\psi_g \bQ_{\cX}|_{X'\setminus\{0\}}=\bQ_{X'\setminus\{0\}}$
and $(\psi_g \bQ_{\cX})_0=\boR\Gamma(F_g,\bC)$.
This means that we have to consider also the contribution of
$H^{d_X+1}(F_g,\bC)$ in order to satisfy the
symmetry and the semicontinuity.
Since the action of the monodromy on $H^{d_X+1}(F_g,\bC)$ is
associated to the function $f$, it is the identity.
So their definition coincides with ours for the
$m_{x,\alpha}$ with $\alpha\notin\bZ$.
However, it seems rather difficult to generalize the construction
in [10] to the case of arbitrary singularities.

\medskip
(v) By the definition of specialization (1.1.2) using the
deformation to the normal cone, we have
$$
\supp\,\Sp_X\bC_Y=N_XY\quad\text{with}\quad\dim N_XY=\dim Y.
$$
The assumption of the next proposition is satisfied in case $X$ has a
0-dimensional embedded component so that
$\dim\,(N_XY)_x=\dim Y$, see also Th.~(4.4) and Ex.~(4.6) below.

\medskip\noindent
{\bf 1.4.~Proposition.} {\it
If $\dim\Lambda=\dim Y$, then $m_{\Lambda,\alpha}=0$ unless
$\alpha\in(0,1]$.
}

\medskip\noindent
{\it Proof.} 
It is enough to show that the restriction of $\Sp_X(\cO_Y,F)$ to a
sufficiently small open subvariety of $\Lambda$ is a variation of
mixed Hodge structure of level $0$ (i.e. the Hodge filtration is
trivial).
Let $n=:\dim\Lambda=\dim Y$, and
$$
(H,F)=i_{\xi}^*\Sp_X(\cO_Y,F)[-n].
$$
This underlies a mixed Hodge structure (where $F^p=F_{-p})$ since
$\xi$ is sufficiently general.
By Remark~(1.3)(iii) we have
$$
\min\{p\,|\,F_pH\ne 0\}\ge 0.
$$
Moreover, the self-duality $\bD(\cO_Y,F)=(\cO_Y,F[n])$ implies
$$
\aligned
\bD(\Sp_X(\cO_Y,F))&=\Sp_X(\cO_Y,F[n]),\\
\bD(H,F)&=(H,F).
\endaligned
$$
Indeed, the second isomorphism follows from the first by using
the duality
$$
i_{\xi}^*\scirc\bD=\bD\scirc i_{\xi}^!
$$
together with
$$
i_{\xi}^!\Sp_X(\cO_Y,F)=i_{\xi}^*\Sp_X(\cO_Y,F[n])[2n].
$$
Here the last isomorphism follows from the fact that
$\Sp_X(\cO_Y,F)[-n]$ underlies a variation of mixed Hodge structure
on a neighborhood of $\xi$ in $\Lambda$.
So the assertion follows.

\bigskip\bigskip
\centerline{\bf 2. Proof of main theorems}

\bigskip
In this section we prove Theorems~1--3.
We first show the following proposition which will be used in the
proof of Theorem~1.

\medskip\noindent
{\bf 2.1.~Proposition.} {\it
Let $(M,F)$ be a filtered $\cD_Y$-module on a smooth variety or
a complex manifold $Y$,
which underlies a mixed Hodge module $\cM$.
Let $p_0=p(M,F)$ in the notation of $(1.3.2)$.
Let $i:X\to Y$ be a closed immersion of a smooth subvariety.
Let $i^*(M,F)$ denote of the pull-back as a complex of filtered
$\cD$-modules. Then

\smallskip\noindent
{\rm (i)} We have $F_{p_0}\cH^ki^*M=0$ for $k\ne 0$, and
$F_{p_0}\cH^0i^*M$ is locally a quotient of
$\cO_X\otimes_{\cO_Y}F_{p_0}M$.

\smallskip\noindent
{\rm (ii)} If $\supp\,F_{p_0}M\subset X$, then
$F_{p_0}\cH^0i^*M=F_{p_0}M$.
}

\medskip\noindent
{\it Proof.}
By the definition of the pull-back,
we may assume that $X$ is defined by $y_1=0$ as in Remark~(1.3)(ii).
Then the assertion (i) for $k\ne 0$ easily follows.
Moreover, the assertion (i) for $k=0$ is reduced to
$$
F_{p_0}M\subset V^0M.
\leqno(2.1.1)
$$
(Indeed, this implies $y_1F_{p_0}M\subset V^{>0}M$, and hence
$F_{p_0}\Gr_V^0M$ is a quotient of $\cO_X\otimes_{\cO_Y}F_{p_0}M$.)
Then (2.1.1) follows from the strict surjectivity of
$$
\partial/\partial y_1:\Gr_V^{\alpha+1}(M,F[1])\to\Gr_V^{\alpha}(M,F)
\quad\text{for}\,\,\,\,\alpha<0,
$$
see [22], 3.2.1.3, where $(F[1])_p=F_{p-1}$ and
$V_{\alpha}=V^{-\alpha}$.

For the assertion (ii) it is enough to show
$F_{p_0}V^{>0}M=0$ under the condition $\supp\,F_{p_0}M\subset X$.
By the definition of the filtration $V$, we have the injectivity of
$$
y_1:V^{\alpha}M\to V^{\alpha+1}M
\quad\text{for}\,\,\,\,\alpha>0.
$$
If $F_{p_0}V^{>0}M\ne 0$, then its support cannot be contained in
$X=\{y_1=0\}$ by the above injectivity of $y_1$.
So the assertion follows.
This completes the proof of Prop.~(2.1).

\medskip\noindent
{\bf 2.2.~Proof of Theorem~1.}
Set $(M,F)=\Sp_X(\cO_Y,F)$.
Applying Prop.~(2.1)(i) to the pull-back by
$$
\ti_x:\{x\}\times\bC^r\to\tY,
$$
we get the assertion~(i), because the pull-back by
$$
i'_{\xi}:\{\xi\}\to\Lambda
$$
is a non-characteristic pull-back (since $\xi\in\Lambda$ is
sufficiently general), and is defined by the pull-back as
$\cO$-modules, see Remark~(1.3)(ii).
(Note that the pull-back by $\Lambda\to\{x\}\times\bC^r$ is
the inverse of the direct image by the closed embedding
if we restrict to a neighborhood of $\xi$ where $\Lambda$ is
smooth.
So this is essentially trivial.)

If $x$ is a sufficiently general point of $Z_{\beta,E}$,
we may assume $Z_{\beta,E}$ is a point by [9], Th.~5.3
taking a sufficiently general member of a family of transversal
slices to $Z_{\beta,E}$.
Here $\dim Y$ and $d_{\beta,E}$ are replaced respectively by
$\dim Y-d_{\beta,E}$ and $0$.
Thus the assertion~(iii) is reduced to (ii),
and the latter follows from Prop.~(2.1)(ii).
This completes the proof of Theorem~1.

\medskip\noindent
{\bf 2.3.~Proof of Theorem~2.}
The first assertion is clear since the restriction of
$\cH^{\ssbull}\Sp_X\bC_Y$ to $S'_j$ underlies a variation of
mixed Hodge structure.

For the second assertion, note that the stratification
$\{S'_j\}$ of $N_XY$  gives a stratification
satisfying Thom's $(a_f)$-condition for the function $t$ in (1.1.2)
since $Y\times\bC^*$ is smooth,
see e.g.\ [9], Prop.~2.17 and the references there.
So the assertions follow from the same argument as in the proofs of
[9], Theorems~1.2 and 5.3.

\medskip\noindent
{\bf 2.4.~Proof of Theorem~3.}
The first assertion on the roots of $b$-function follows from [4],
Cor.~2.8, and the second assertion on the jumping coefficients
follows from Theorem~1(i).
Indeed, if the second assertion does not hold, then the degree $i$
part of $\cM(\beta,x)$ vanishes for $i\gg 0$ (since the graded
algebra $\ocA$ is generated by the degree 1 part).
But this implies that $\supp\,\cM(\beta,x)\subset\{0\}$ in
$\Spec\,\ocA(x)=(T_XY)_x$, and it contradicts Theorem~1(i).
This completes the proof of Theorem~3.

\medskip\noindent
{\bf 2.5.~Remark.} Let
$E_{\Lambda}=\{\alpha\,|\,m_{\Lambda,\alpha}\ne 0\}$, and
$R_{f,x}$ be the set of the roots of $b_{f,x}(-s)$, where
$b_{f,s}(s)$ is defined for the germ $(X,x)$, see [4].
Then the first assertion of Theorem~3 is equivalent to
$$
\hbox{$\bigcup_{\Lambda}\exp(-2\pi iE_{\Lambda})\subset
\exp(-2\pi iR_{f,x})$},
$$
where $\Lambda$ runs over the irreducible components of $(N_XY)_x$.
However, the equality does not always hold (e.g. if $f=x^2y$)
unless we take the union over the irreducible components
$\Lambda$ of $(N_XY)_y$ for any $y\in Z$ sufficiently near $x$.

\bigskip\bigskip
\centerline{\bf 3. Product formula}

\bigskip
In this section we show a product formula which implies the
well-definedness of the spectrum.

\medskip\noindent
{\bf 3.1.~Cartesian product.}
For $a=1,2$, let $Y_a$ be a smooth complex algebraic variety or
a complex manifold, and $X_a$ be a closed subvariety of $Y_a$.
Let $X=X_1\times X_2$, $Y=Y_1\times Y_2$ with the projection
$pr_a$ to the $a$-th factor.
Let $\cI_a$ denote the ideal of $X_a\subset Y_a$.
Then $\cI_X=pr_1^*\cI_1+pr_2^*\cI_2$, and
$$
\cI_X^i=\msum_{p+q=i}\,\cI_1^p\boxtimes\cI^q=
\msum_{p+q=i}\,pr_1^*\cI_1^p\cap pr_2^*\cI_2^q,
\leqno(3.1.1)
$$
i.e. the filtration $\{\cI_X^i\}_{i\in\bN}$ is the convolution of
the filtrations $\{pr_1^*\cI_1^p\}_{p\in\bN}$ and
$\{pr_2^*\cI_2^q\}_{q\in\bN}$.
Here the last isomorphism of (3.1.1) is shown by taking the exterior
product of the exact sequences
$$
\aligned
0\to\cI_1^p\to\cO_{Y_1}\to\cO_{Y_1}/\cI_1^p\to 0,\\
0\to\cI_2^q\to\cO_{Y_2}\to\cO_{Y_2}/\cI_2^q\to 0,
\endaligned
$$
which gives a diagram of short exact sequences by the exactness of
exterior product.
By a similar argument, we get then (see also [1], 3.1.2.9)
$$
\cI_X^i/\cI_X^{i+1}=
\mopl_{p+q=i}\,(\cI_1^p/\cI_1^{p+1}\boxtimes\cI_2^q/\cI_2^{q+1}),
\leqno(3.1.2)
$$
i.e.
$$
\mopl_{i\in\bN}\,\cI_X^i/\cI_X^{i+1}=
\big(\mopl_{p\in\bN}\,\cI_1^p/\cI_1^{p+1}\big)\boxtimes
\big(\mopl_{q\in\bN}\,\cI_2^q/\cI_2^{q+1}\big),
$$
and hence
$$
N_XY=N_{X_1}Y_1\times N_{X_2}Y_2.
\leqno(3.1.3)
$$

\medskip\noindent
{\bf 3.2.~Proposition.} {\it
With the above notation we have a canonical isomorphism
$$
\Sp_X(\cO_Y,F)=
\Sp_{X_1}(\cO_{Y_1},F)\boxtimes\Sp_{X_2}(\cO_{Y_2},F),
\leqno(3.2.1)
$$
where the action of $T_s$ on the left-hand side corresponds
to $T_s\boxtimes T_s$ on the right-side hand.
}

\medskip\noindent
{\it Proof.}
Let $M_a$ be the direct image of $\cO_{Y_a}$ by the graph
embedding as in (1.1).
Let
$$
(M,F)=(M_1,F)\boxtimes(M_2,F).
$$
These have the filtration $V$ of Kashiwara [14] and Malgrange [18].
Using the same argument as in the proof of (3.1.1), we see that
the filtration $V$ on $M$ is the convolution of $pr_1^*V$ and
$pr_2^*V$, i.e.
$$
V^{\alpha}M=\msum_{\alpha_1+\alpha_2=\alpha}\,
V^{\alpha_1}M_1\boxtimes V^{\alpha_2}M_2.
$$
(Indeed, the filtration defined by the right-hand side satisfies
the conditions of the filtration $V$.)
Then we have
$$
\Gr_V^{\alpha}(M,F)=\mopl_{\alpha_1+\alpha_2=\alpha}\,
\Gr_V^{\alpha_1}(M_1,F)\boxtimes\Gr_V^{\alpha_2}(M_2,F).
$$
So the assertion follows.

\medskip\noindent
{\bf 3.3.~Corollary.} {\it
With the above notation, write
$\hSp(X_a,\Lambda_a)=\msum_{\alpha}m_{a,\Lambda_a,\alpha}t^{\alpha}$
for an irreducible component $\Lambda_a$ of $(N_{X_a}Y_a)_{x_a}$,
and $\hSp(X,\Lambda)=\msum_{\alpha}m_{\Lambda,\alpha}
t^{\alpha}$ for $\Lambda=\Lambda_1\times\Lambda_2$ under the
isomorphism $(3.1.3)$.
Then
$$
m_{\Lambda,\alpha}=\msum_{\alpha_1\mp\alpha_2=\alpha}\,
m_{1,\Lambda_1,\alpha_1}m_{2,\Lambda_2,\alpha_2}.
\leqno(3.3.1)
$$
where $\alpha_1\mp\alpha_2$ is defined to be
$\alpha_1+\alpha_2$ if
$\lceil\alpha_1\rceil-\alpha_1+\lceil\alpha_2\rceil-\alpha_2\ge 1$,
and $\alpha_1+\alpha_2-1$ otherwise.
Here $\lceil\alpha\rceil$ is the smallest integer which is greater
than or equal to $\alpha$.
}

\medskip\noindent
{\it Proof.}
Let $\xi=(\xi_1,\xi_2)\in N_XY=N_{X_1}Y_1\times N_{X_2}Y_2$.
By Prop.~(3.2) we have
$$
i^*_{\xi}\Sp_X(\cO_Y,F)=
i^*_{\xi_1}\Sp_{X_1}(\cO_{Y_1},F)\boxtimes
i^*_{\xi_2}\Sp_{X_2}(\cO_{Y_2},F).
$$
So the assertion follows.

\medskip\noindent
{\bf 3.4.~Corollary.} {\it
The spectrum $\hSp(X,\Lambda)$ is essentially independent of $Y$
using the isomorphism $(3.1.3)$ for $X_2=pt$.
}

\medskip\noindent
{\it Proof.} Since the spectrum is defined analytically locally,
it is enough to compare the embedding $X\to Y$ with
$$
X\to Y=Y\times\{0\}\to Y\times\bC.
$$
So the assertion follows from (3.3), since
$\hSp(X_2,\Lambda_2)=t$ in the case $X_2$ is a (reduced) point,
$\dim Y_2=1$, and $\Lambda_2=N_{X_2}Y_2$.

\bigskip\bigskip
\centerline{\bf 4. Monomial ideal case}

\bigskip
In this section, we treat the monomial ideal case, and prove
Th.~(4.4).

\medskip\noindent
{\bf 4.1.~Notation.}
Assume $Y=\bC^n$, and $X$ is defined by a monomial ideal $\fa$ of
$\bC[x]=\bC[x_1,\dots, x_n] $.
We denote by $x^{\nu}$ the monomial corresponding to ${\nu}\in\bN^n$.
Let $\Gamma_{\fa}\subset\bN^n$ be the semigroup corresponding
to $\fa $, i.e.
$$
\Gamma_{\fa} = \{{\nu}\in\bN^n\mid x^{\nu}\in\fa\}.
$$
Let $P_{\fa}$ be the convex hull of $\Gamma_{\fa}$ in $\bR_{\ge 0}^{n}$
which is called the Newton polyhedron of $\fa $.
Let $J(Y,\alpha X)\subset \bC[x]$ denote the multiplier ideals of $\fa$.
Set ${\mathbf 1}=(1,\dots,1)$, and
$$
U(\alpha)=\{\nu\in\bN^n\,|\,{\nu}+{\mathbf 1}\in(\alpha+\varepsilon)
P_{\fa}\,\,\,\text{with}\,\,\,0<\varepsilon\ll 1\}.
$$
Then we have by Howald~[12]
$$
J(Y,\alpha X)=\msum_{\nu\in U(\alpha)} \bC\,x^{\nu}.
\leqno(4.1.1)
$$
See [5] for a combinatorial description of the roots of the
$b$-function.

For a $(n-1)$-dimensional compact face $\sigma$ of $P_{\fa}$,
let $L_{\sigma}$ be the linear function such that $L_{\sigma}^{-1}(1)
\supset \sigma$. Let $c_{\sigma}$ be the smallest positive integer
such that $c_{\sigma}L_{\sigma}$ has integral coefficients.
Let
$$
G'_{\sigma}=\bZ^n\cap L^{-1}_{\sigma}(0),
$$
and $G_{\sigma}$ be the subgroup generated by $\nu-\nu'$ for
$\nu,\nu'\in\Gamma_{\fa}\cap\sigma$.
Set
$$
e_{\sigma}=|G'_{\sigma}/G_{\sigma}|.
$$

For a face $\sigma$ of $P_{\fa}$, let
$\oB_{\sigma}\subset\bC[x]$ be the $\bC$-subalgebra generated by
$x^{\nu}$ for
$\nu\in\sigma\cap\Gamma_{\fa}$. Let
$$
\oB=\msum_{\sigma}\,\oB_{\sigma}\subset\bC[x],
$$
where the multiplication of $x^{\nu}$ and $x^{\nu'}$ in $\oB$ is
given by $x^{\nu+\nu'}$ if $x^\nu,x^{\nu'}\in \oB_{\sigma}$ for some
$\sigma$, and it vanishes otherwise. Set
$$
\oA=\mopl_{i\in\bN}\,\fa^i/\fa^{i+1}.
$$

\medskip
With the above notation, we have the following:

\medskip\noindent
{\bf 4.2.~Proposition.}
$\oA_{\red}=\oB$.

\medskip\noindent
{\it Proof.}
For $\nu\in\bN^n$, set
$$
v(\nu)=\min\{L_{\sigma}(\nu)\},
$$
where $\sigma$ runs over the $(n-1)$-dimensional faces of
$P_{\fa}$.
Let $\hsig$ be the cone generated by $\sigma$ in the real vector space
$\bR^n$.
Then
$$
\bR_{>0}^n\subset\mcup_{\sigma}\hsig\subset\bR_{\ge 0}^n,
$$
and
$$
v(\nu)=\begin{cases}
L_{\sigma}(\nu) &\text{if}\,\,\, \nu\in\hsig,\\
0 &\text{if}\,\,\, \nu\notin\mcup_{\sigma}\hsig.
\end{cases}
\leqno(4.2.1)
$$
So we get for $\nu\in\bR_{>0}^n$
$$
\nu\in \alpha P_{\fa}\Leftrightarrow v(\nu)\ge\alpha.
\leqno(4.2.2)
$$
Note that
$$
v(\nu+\nu')=L_{\sigma''}(\nu+\nu')\ge v(\nu)+v(\nu')\quad
\hbox{if}\,\,\, \nu+\nu'\in\hsig''.
\leqno(4.2.3)
$$

Let $C$ be a positive number such that for any $\nu\in\bN^n$ and
$k\in\bN$
$$
x^{\nu}\in \fa^k\,\,\,\text{if}\,\,\,v(\nu)\ge k+C.
\leqno(4.2.4)
$$
For the existence of such $C$, it is enough to show
$$
\msum^k\Gamma_{\fa}+\bR_{\ge 0}^n\supset(k+C)P_{\fa},
\leqno(4.2.5)
$$
where
$\msum^kS:=\{\sum_{i=1}^k v_i\,|\,v_i\in S\}$ for $S\subset \bR^n$.
Then the assertion is reduced to
$$
\msum^k(\sigma\cap\Gamma_{\fa})+\bR_{\ge 0}^n\supset
\hsig\cap(k+C)P_{\fa}\quad\hbox{for any}\,\,\,\sigma.
\leqno(4.2.6)
$$
If $\sigma$ is not compact, it is the union of
$\sigma'+\bR_{\ge 0}^I$ for compact faces $\sigma'$ of $\sigma$,
where $I$ is the subset of $\{1,\dots,n\}$ such that $\sigma$
is stable by adding the $i$th unit vector for $i\in I$.
So the assertion is reduced to the case $\sigma$ compact.
By increasing induction on $k$, it is further is reduced to
$$
(\sigma\cap\Gamma_{\fa})+\hsig\supset\hsig\cap\alpha P_{\fa}\quad
\hbox{if}\,\,\,\alpha\gg 0,
\leqno(4.2.7)
$$
since $v(\nu)=1$ for $\nu\in\sigma\cap\Gamma_{\fa}$.
Then the assertion is proved by replacing $\sigma$ with a simplex
$\sigma'$ defined by vertices of $\sigma$.

Set
$$
N=[(C+1)/\varepsilon]+1\quad\text{with}\quad
\varepsilon=\min\{c_{\sigma}^{-1}\},
$$
where $\sigma$ runs over the $(n-1)$-dimensional faces
of $P_{\fa}$.
Note that
$$
v(\nu)\ge k+\varepsilon \quad\hbox{if}\quad
v(\nu)>k \,\,\,\hbox{with}\,\,\,
\nu\in\bN^n.
\leqno(4.2.8)
$$
Let $I_k$ be the ideal of $\bC[x]$ generated by $x^{\nu}\in\fa^k$
with $v(\nu)>k$.
Set
$$
\oI=\mopl_{k\in\bN}\,I_k/\fa^{k+1}.
$$
Then $\oI$ is an ideal of $\oA$ such that $\oI{}^N=0$ and
$\oA/\oI=\oB$ by the above arguments.
So the assertion follows.

\medskip\noindent
{\bf 4.3.~Corollary.} {\it
{\rm (i)}
The irreducible components $\Lambda$ of $N_XY$ are given by
$\Spec\,\oB_{\sigma}$ for $(n-1)$-dimensional faces $\sigma$ of
$P_{\fa}$,
and we have for any faces $\sigma,\sigma'$ of $P_{\fa}$
$$
\Spec\,\oB_{\sigma}\cap\Spec\,\oB_{\sigma'}=
\Spec\,\oB_{\sigma\cap\sigma'}.
$$
{\rm (ii)} The irreducible components $\Lambda$ of $(N_XY)_0$
are given by $\Spec\,\oB_{\sigma}$ for $(n-1)$-dimensional compact
faces $\sigma$ of $P_{\fa}$.
}

\medskip\noindent
{\it Proof.}
The first assertion of (i) is clear by Prop.~(4.2). Let
$$
J_{\sigma}=\Ker(\oB\to\oB_{\sigma}).
$$
This is generated over $\bC$ by $x^{\nu}$ for
$\nu\in\bN^n$ such that $x^{\nu}\in\oB$ and $\nu\notin\hsig$.
So we get
$$
J_{\sigma\cap\sigma'}=J_{\sigma}+J_{\sigma'},
$$
and the last assertion of (i) follows.
For the assertion (ii), note that the image of $\Spec\,\oB_{\sigma}$
in $Y$ is a point if and only if $\sigma$ is compact.
Then the assertion (ii) follows from (i).

\medskip\noindent
{\bf 4.4.~Theorem.} {\it
For $\Lambda=\Spec\,\oB_{\sigma}$ with $\sigma$ compact,
we have in the notation of {\rm (4.1)}
}
$$
\hSp(X,\Lambda)=\msum_{i=1}^{c_{\sigma}}\,e_{\sigma}t^{i/c_{\sigma}}.
$$

\medskip\noindent
{\it Proof.}
Let $\oB'_{\sigma}$ be the localization of $\oB_{\sigma}$ by the
monomials $x^{\nu}$ in $\oB_{\sigma}$.
Then
$$
\rank_{\oB'_{\sigma}}\cM(\beta)\otimes_{\oA}\oB'_{\sigma}=
\begin{cases}
e_{\sigma} &\text{if}\,\,\, \beta\in c_{\sigma}^{-1}\bZ,\\
0  &\text{otherwise}.
\end{cases}
\leqno(4.4.1)
$$
So the assertion follows from Theorems~1 together with
Prop.~(1.4), (4.2) and Cor.~(4.3).

\medskip\noindent
{\bf 4.5.~Comparison.}
With the notation of (4.1), let $f=\sum_{\nu}c_{\nu}x^{\nu}\in
\bC\{\!\{x\}\!\}$.
Assume $f$ is non-degenerate with respect to the Newton boundary
$\partial P_{\fa}$, see [15].
Let $D=f^{-1}(0)\subset (Y,0)=(\bC^n,0)$.
Assume
$$
X_{\red}=\{0\}\quad \text{so that}\quad \Sing D=\{0\}.
$$
Then we may assume $f\in\bC[x]$ by the finite determination property
as is well-known.
Let $\JC(Y,X)$, $\JC(Y,D)$ denote the sets of jumping coefficients.
With the notation in the proof of Prop.~(4.2) we have by Howald [12]
(see (4.1.1) above)
$$
\JC(Y,X)=\big\{v(\nu)\,\big|\,\nu\in\bZ_{>0}^n\big\}.
\leqno(4.5.1)
$$
Note that $\JC(Y,D)$ contains always 1, and we have the periodicity
$$
\JC(Y,D)=\big(\JC(Y,D)\cap(0,1]\big)+\bN.
\leqno(4.5.2)
$$

Write $\Sp(f,0)=\sum_{\alpha}m_{f,\alpha}t^{\alpha}$ (the spectrum of
$f$).
Define the set of exponents by
$$
E(D,0)=\{\alpha\,|\, m_{f,\alpha}\ne 0\},\quad
E(X,\Lambda)=\{\alpha\,|\, m_{\Lambda,\alpha}\ne 0\}.
$$
Then we have

$$
\begin{array}{ccc}
\JC(Y,D)\cap(0,1) & \overset{(1)}{=\!=} &
\JC(Y,X)\cap(0,1)\\
{\scriptstyle (2)\big|\!\big|} &  &
\cap{\scriptstyle (3)}\\
E(D,0)\cap(0,1) & \underset{(4)}{\subset} &
\mcup_{\Lambda}E(X,\Lambda)\cap(0,1)\\
\end{array}
\leqno(4.5.3)
$$
Indeed, we have the equality (1) by Howald [13], and
(2) by Budur [3].
The inclusion (3) follows from Theorem~1 and (4.4.1)
or Theorem~(4.4) and (4.5.1).
The spectrum of $f$ is calculated by Steenbrink [27] (see also
[21], [30]). This implies that the composition of (1) and (2)
is the equality, since the exponents at most 1 are given by
restricting the right-hand side of (4.5.1).
Combined with Theorem~(4.4), the last assertion also implies (4).
Note that equality does not necessarily hold for (3) and (4),
see Ex.~(4.6--7) below.

As for the $b$-function, it is well known that the roots of the
$b$-function are not determined by the Newton boundary even if
$f$ is non-degenerate in the sense of [15].
For example, the roots of the $b$-function of
$f_{\lambda}=x^5+y^4+\lambda x^3y^2$ depend on whether
$\lambda=0$ or not.
Indeed, for $\lambda\ne 0$, $\frac{11}{20}$ is a root of
$b_{f_{\lambda}}(-s)$ and $\frac{31}{20}$ is not.
For $\lambda=0$, see (4.7.1) below.

\medskip\noindent
{\bf 4.6.~Example.}
Assume $Y=\bC^n$ and $f=(f_1,\dots,f_n)$ with $f_i=x_i^{m_i}$,
where $x_1,\dots,x_n$ are the coordinates of $\bC^n$ and
the $m_i$ are positive integers.
Let $\sigma$ be the unique $(n-1)$-dimensional compact face of
$P_{\fa}$,
i.e. the convex full of
$$
\{m_1\be_1,\dots,m_n\be_n\},
$$
where $\be_i$ is the $i$-th unit vector of $\bR^n$.
Then $L_{\sigma}=\msum_i\,x_i/m_i$, and
$$
c_\sigma=\hbox{LCM}(m_1,\dots,m_n),\quad
e_{\sigma}=m_1\cdots m_n/c_{\sigma}.
\leqno(4.6.1)
$$
Indeed, let $\oG_{\sigma}$ and $\oG'_{\sigma}$ denote respectively
the injective image of $G_{\sigma}$ and $G'_{\sigma}$ by the
projection $\bZ^n\to\bZ^{n-1}$ to the first $n-1$ factors.
Set
$$
b=\hbox{LCM}(m_1,\dots,m_{n-1}),\,\,\,
d=\hbox{GCD}(b,m_n).
$$
Then $b=b'd$, $m_n=m'_nd$ with $b',m'_n\in\bN$.
This implies $c_{\sigma}=m_nb'$, and
$$
|\bZ^{n-1}/\oG_{\sigma}|=m_1\cdots m_{n-1},\quad
|\bZ^{n-1}/\oG'_{\sigma}|=b'.
$$
So the assertion for $e_{\sigma}$ follows.
The assertion for $c_{\sigma}$ is clear.

Set
$$
E=\{(a_1,\dots,a_n)\in\bN^n\,|\,a_i\in[1,m_i]\}.
$$
Then
$$
\aligned
\JC(Y,X)&=\Big\{\sum_{i=1}^n\frac{a_i}{m_i}\,\Big|\,
a_i\in\bZ_{>0}\Big\},\\
\hSp(X,0)&=\sum_{i=1}^{c_{\sigma}}e_{\sigma}t^{i/c_{\sigma}},\\
b_f(s)&=\Big[\prod_{(a_1,\dots,a_n)\in E}\Big(s+\sum_{i=1}^n
\frac{a_i}{m_i}\Big)\Big]_{\red}.
\endaligned
\leqno(4.6.2)
$$
Here $[\prod_j(s+\beta_j)^{n_j}]_{\red}=\prod_j(s+\beta_j)$
if the $\beta_j$ are mutually different and $n_j\in\bZ_{>0}$.
The assertion on the spectrum follows from Th.~(4.4) or
Cor.~(3.3) using (4.6.1).
The other assertions follow from (4.5.1) and [4], Th.~5.
This example shows that it is not necessarily easy to determine
$i$ and $j_0$ in Theorem~3 in general.

\medskip\noindent
{\bf 4.7.~Example.}
Let $f=\sum_i x_i^{m_i}$ and $D=f^{-1}(0)\subset Y=\bC^n$.
Set
$$
\tE=\{(a_1,\dots,a_n)\in\bN^n\,|\,a_i\in[1,m_i-1]\}.
$$
Then
$$
\aligned
\JC(Y,D)\cap(0,1]&=\Big\{\sum_{i=1}^n\frac{a_i}{m_i}
\,\Big|\,a_i\in\bZ_{>0}\Big\}\cap(0,1],\\
\Sp(D,0)&=\prod_{i=1}^n(t-t^{1/m_i})/(t^{1/m_i}-1),\\
\tb_f(s)&=\Big[\prod_{(a_1,\dots,a_n)\in \tE}\Big(s+\sum_{i=1}^n
\frac{a_i}{m_i}\Big)\Big]_{\red},
\endaligned
\leqno(4.7.1)
$$
where $\tb_f(s)=b_f(s)/(s+1)$.
The assertions on $\JC(Y,D)$ and $\Sp(D,0)$ follow from [13] and
[26].
The assertion on $\tb_f(s)$ is an unpublished result
of Kashiwara asserting that in the isolated weighted homogeneous
singularity case, the roots of $\tb_f(-s)$ coincide with the
exponents and have multiplicity 1. (This also follows from [17].)
Note that the formula for $\Sp(D,0)$ holds in the isolated
weighted homogeneous singularity case if we replace $1/m_i$ by
the weights $w_i$.

\medskip\noindent
{\bf 4.8.~Remark.}
In the monomial ideal case, $j_0$ in Theorem~3 is bounded by
$n-1$.
Indeed, with the notation of the proof of Prop.~(4.2) we have for
$\beta\in\bQ\cap(0,1]$ and a face $\sigma$ of $P_{\fa}$
$$
\min\{L_{\sigma}(\nu)\,|\,\nu\in\bZ^n\cap\hsig,\,\,
L_{\sigma}(\nu)-\beta\in\bZ\}\le n.
\leqno(4.8.1)
$$
To show this, we may replace $\hsig$ with
$\hsig'+\sum_{i\in I}\bR_{\ge 0}^I$
where $\sigma'$ is a simplex defined by vertices $\{v_i\}$
of a compact face of $\sigma$, and $I$ is as in the proof of
Prop.~(4.2).
Set
$$
D_{\sigma'}=\Int\hsig'\setminus\mcup_{i}(\Int\hsig'+v_i),
$$
where $\Int\hsig'$ is the interior of $\hsig'$.
Then
$$
\Int\hsig'=\mcup_{\nu\in\bN}(D_{\sigma}+\msum_i \nu_iv_i),
$$
and (4.8.1) follows.
By a similar argument, $\JC(Y,X)$ is stable by adding any
positive integers in the monomial ideal case.
Note that $j_0=n-1$ if the $m_i$ in (4.6) are mutually prime.
In general it is unclear whether $j_0$ is always bounded by $n-1$.

\bigskip
\ver


\begin{thebibliography}{99}

\bibitem{[1]}
Beilinson, A., Bernstein, J. and Deligne, P., Faisceaux pervers,
Ast\'erisque 100, Soc. Math. France, Paris, 1982.

\bibitem{[2]}
Borel, A. et al., Algebraic $\cD$-modules,
Perspectives in Math. 2, Academic Press, 1987.

\bibitem{[3]}
Budur, N., On Hodge spectrum and multiplier ideals,
Math. Ann. 327 (2003), 257--270.

\bibitem{[4]}
Budur, N., Mustata, M. and Saito, M.,
Bernstein-Sato polynomials of arbitrary varieties,
Compos. Math. 142  (2006), 779--797.

\bibitem{[5]}
Budur, N., Mustata, M. and Saito, M.,
Combinatorial description of the roots of the Bernstein-Sato
polynomials for monomial ideals, Comm. Algebra 34 (2006), 4103--4117.

\bibitem{[6]}
Budur, N. and Saito, M., Multiplier ideals,
$ V $-filtration, and spectrum,
J. Algebraic Geom. 14 (2005), 269--282.

\bibitem{[7]}
Deligne, P., Th\'eorie de Hodge I, Actes Congr\`es Intern.
Math., Part 1 (1970), 425--430; II, Publ. Math. IHES, 40 (1971),
5--58; III, ibid. 44 (1974), 5--77.

\bibitem{[8]}
Deligne, P., Le formalisme des cycles \'evanescents, in SGA7 XIII
and XIV, Lect. Notes in Math. 340, Springer, Berlin, 1973,
pp. 82--115 and 116--164.

\bibitem{[9]}
Dimca, A., Maisonobe, Ph., Saito, M., and Torrelli, T.,
Multiplier ideals, $V$-filtrations and transversal sections,
Math. Ann. 336 (2006), no. 4, 901--924.

\bibitem{[10]}
Ebeling, W. and Steenbrink, J.H.M.,
Spectral pairs for isolated complete intersection singularities,
J. Algebraic Geom. 7 (1998), 55--76. 

\bibitem{[11]}
Gyoja, A.,
Bernstein-Sato's polynomial for several analytic functions,
J. Math. Kyoto Univ. 33 (1993), 399--411.

\bibitem{[12]}
Howald, J., Multiplier ideals of monomial ideals, Trans. Amer.
Math. Soc. 353 (2001), 2665--2671.

\bibitem{[13]}
Howald, J., Multiplier ideals of sufficiently general polynomials
(math.AG/0303203).

\bibitem{[14]}
Kashiwara, M., Vanishing cycle sheaves and holonomic systems of
differential equations,
Lect. Notes in Math. 1016, Springer, Berlin, 1983, pp. 134--142.

\bibitem{[15]}
Kouchinirenko, A., Poly\`edres de Newton et nombres de Milnor,
Inv. Math. 32 (1976), 1--31.

\bibitem{[16]}
Lazarsfeld, R.,
Positivity in algebraic geometry II,
Springer-Verlag, Berlin, 2004.

\bibitem{[17]}
Malgrange, B., Le polyn\^ome de Bernstein d'une
singularit\'e isol\'ee, in Lect. Notes in Math. 459, Springer,
Berlin, 1975, pp. 98--119.

\bibitem{[18]}
Malgrange, B., Polyn\^ome de Bernstein-Sato et cohomologie
\'evanescente, Analysis and topology on singular spaces, II, III
(Luminy, 1981), Ast\'erisque 101--102 (1983), 243--267.

\bibitem{[19]}
Parameswaran, A.J., Monodromy fibration of an isolated complete
intersection singularity, Proc. Indo-French Conference on Geometry
(Bombay, 1989),  Hindustan Book Agency, Delhi, 1993, pp. 123--134.

\bibitem{[20]}
Sabbah, C.,
Proximit\'e \'evanescente, I. La structure polaire d'un
$ D $-module, Compos. Math. 62 (1987), 283--328;
II. Equations fonctionnelles pour plusieurs fonctions analytiques,
ibid. 64 (1987), 213--241.

\bibitem{[21]}
Saito, M., Exponents and Newton polyhedra of isolated hypersurface
singularities, Math. Ann. 281 (1988), 411--417.

\bibitem{[22]}
Saito, M., Modules de Hodge polarisables, Publ. RIMS, Kyoto
Univ. 24 (1988), 849--995.

\bibitem{[23]}
Saito, M., Mixed Hodge modules, Publ. RIMS, Kyoto Univ. 26
(1990), 221--333.

\bibitem{[24]}
Saito, M., On $b$-function, spectrum and rational singularity,
Math. Ann. 295 (1993), 51--74.

\bibitem{[25]}
Saito, M., Multiplier ideals, $b$-function, and spectrum of a
hypersurface singularity, preprint (math.AG/0402363).

\bibitem{[26]}
Steenbrink, J.H.M.,
Intersection form for quasi-homogeneous singularities,
Compos. Math. 34 (1977), 211--223.

\bibitem{[27]}
Steenbrink, J.H.M., Mixed Hodge structure on the vanishing
cohomology, in Real and Complex Singularities (Proc. Nordic
Summer School, Oslo, 1976) Alphen a/d Rijn: Sijthoff \& Noordhoff
1977, pp. 525--563.

\bibitem{[28]}
Steenbrink, J.H.M., The spectrum of hypersurface singularity,
Ast\'erisque 179--180 (1989), 163--184.

\bibitem{[29]}
Steenbrink, J.H.M., Spectra of $\mathcal K$-unimodal isolated
singularities of complete intersections,
Singularity theory (Liverpool, 1996), London Math. Soc.
Lecture Note Ser., 263, Cambridge Univ. Press, Cambridge, 1999,
pp. 151--162.

\bibitem{[30]}
Varchenko, A.N, Khovanski, A.G.,
Asymptotic behavior of integrals over vanishing cycles and
the Newton polyhedron,
Soviet Math. Dokl. 32 (1985), 122--127.

\bibitem{[31]}
Verdier, J.-L., Sp\'ecialisation de faisceaux et monodromie
mod\'er\'ee, Analysis and topology on singular spaces, II,
III (Luminy, 1981), Ast\'erisque 101--102 (1983), 332--364.

\end{thebibliography}
\end{document}